\documentstyle[12pt]{article}

\def\beq{\begin{equation}}
\def\ee{\end{equation}}
\def\bib#1{[{\ref{#1}}]}
\def\psiba{\overline{\psi}}
\def\vphiba{\overline{\varphi}}
\def\zba{\overline{z}}
\def\pba{\overline{p}}

\def\CC{\mbox{\boldmath $C$}}
\def\RR{R}
\def\SS{S}

\newtheorem{proposition}{Proposition}
\newtheorem{corollary}{Corollary}
\newtheorem{conjecture}{Conjecture}
\begin{document}           

            \title{
Weierstrass representations for \\
surfaces in 4D spaces and their \\
integrable deformations via DS hierarchy}
 \author{B.G. Konopelchenko \\
{\em Dipartimento di Fisica, Universit\'{a} di Lecce, 
73100 Lecce, Italy} } 

 \maketitle

      \begin{abstract}

Generalized Weierstrass representations for 
generic surfaces conformally immersed 
into four-dimensional Euclidean and pseudo-Euclidean spaces of different 
signatures are presented. Integrable deformations of surfaces in these 
spaces generated by the Davey-Stewartson hierarchy of integrable equations 
are proposed. Willmore functional of a surface is invariant under such 
deformations. 
      \end{abstract}

\hfill




\newpage

\section{Introduction}
\setcounter{equation}{0}

Surfaces and their deformations (dynamics) were the subject 
of intensive study for a long time
both in mathematics and physics.
Theories of immersion and deformations of surfaces have been developed 
in the classical differential geometry since the
last century (see {\em e.g.} 
\bib{r1}-\bib{r3}). They continue to be an important part of the
differential geometry (see \bib{r4}-\bib{r5}). 
In physics, surfaces are key ingredients in a broad variety of phenomena 
from surface waves in hydrodynamics to world-sheets in string theory
(see {\em e.g.} \bib{r6}-\bib{r8}). 

Analytic methods to study surfaces and their properties are of great 
interest both in mathematics and in physics. A classical example of such an 
approach is given by the Weierstrass representation for minimal surfaces
(see {\em e.g.} \bib{r1}-\bib{r3}). This representation allows us to 
construct any minimal surface in the three-dimensional Euclidean space 
$\RR^3$ via two holomorphic functions. It is the most powerful tool for an 
analysis of minimal surfaces. 

Extensions of the old Weierstrass representation to generic nonminimal 
surfaces in $\RR^3$ have been given in \bib{r9} and \bib{r10}-\bib{r11}. 
Looking differently these two extensions have occurred to be equivalent. The 
generalized Weierstrass representation proposed in \bib{r10}-\bib{r11} 
starts with the linear system 
\beq 
\begin{array}{l}
\psi_{z}=p \varphi \;\; ,\\
\varphi_{\zba}=-p \psi
\end{array}
\label{eq:1.1}
\ee
where $z=x+iy$, bar means complex conjugation, $\psi(z,\zba)$, 
$\varphi(z,\zba)$ are complex-valued functions while $p(z,\zba)$ is a 
real-valued one. Then conformal immersion of a surface into $\RR^3$ with 
coordinates $X^1$, $X^2$, $X^3$ is defined by the formulae 
\bib{r10}-\bib{r11}

\begin{eqnarray}
X^2 & + & i X^1  = i \int_{\Gamma} {\left( \overline{\psi}^2 dz'-
\overline{\varphi}^2 d\zba' \right)} \; ,  \nonumber \\
X^3 &  = &  
\int_{\Gamma} { \left( \psiba \varphi  dz' + \vphiba \psi  d\zba'
\right) } \;\; 
\label{eq:1.2}
\end{eqnarray}
where $\Gamma$ is a contour in $\CC$. The induced metric on a surface is 
given by
\beq
ds^2=\left( |\psi|^2+|\varphi|^2 \right)^2 dz d\zba ,
\label{eq:1.3}
\ee
the mean curvature is
\beq
H=2\frac{p}{|\psi|^2+|\varphi|^2} \; \; .
\label{eq:1.4}
\ee
and the Willmore functional (see {\em e.g.} \bib{r12})
$W=\int {H^2 \left[ dS \right] }$ is equal to
\beq
W=4 \int{ p^2 dx \; dy}
\label{eq:1.5}
\ee

An advantage of the generalized Weierstrass formulae 
(\ref{eq:1.1}), (\ref{eq:1.2}) 
is also that it allows to construct new class of deformations 
of surfaces in $\RR^3$ \bib{r10}-\bib{r11}. 
They are generated by evolution of $p$ via the modified 
Veselov-Novikov equation and of
$\psi$, $\varphi$ via certain linear equations
. A characteristic feature of these integrable deformations is that the 
Willmore functional remains invariant.

The generalized 
Weierstrass representation (\ref{eq:1.1}), (\ref{eq:1.2}) 
has been proved to be an 
effective tool to study generic surfaces in $\RR^3$ and their 
deformations. In differential geometry its use has allowed to obtain 
several interesting results both of local and global character
(see {\em e.g.} \bib{r13}-\bib{r20}). In physics, it has been applied to 
study of various problems in theory of liquid membranes, $2D$ gravity and 
string theory \bib{r13}, \bib{r21}-\bib{r23}.
In the string theory the functional 
$W=\int {H^2 \left[ dS \right] }$ is known as the Polyakov 
extrinsic action and in membrane theory it is the Helfrich free energy
\bib{r6}-\bib{r8}).

An extension of the Weierstrass representation to multidimensional spaces 
would be of a great interest. In physics, a strong motivation lies in the 
Polyakov string integral over surfaces in multidimensional spaces 
\bib{r6}-\bib{r8}). Theory of immersion of surfaces
into four-dimensional spaces is an important part of the contemporary 
differential geometry too (see {\em e.g.} \bib{r4}, \bib{r5}, 
\bib{r12}, \bib{r24}-\bib{r27}). 

In this paper we present extensions of the generalized Weierstrass 
representation (\ref{eq:1.1}), (\ref{eq:1.2}) to the cases of generic 
surfaces immersed into the four-dimensional spaces
$\RR^4$, $\RR^{3,1}$ and $\RR^{2,2}$ with the metrics 
$g_{ik}=diag(1,1,1,1)$, $g_{ik}=diag(1,1,1,-1)$ and
$g_{ik}=diag(1,1,-1,-1)$, respectively. A basic linear system consists of 
the two two-dimensional Dirac equations while the formulae of immersion are 
of the type (\ref{eq:1.2}). The induced metric, mean curvature and Willmore 
functional are given by formulae similar to (\ref{eq:1.3})-(\ref{eq:1.5}).

Integrable deformations of surfaces are generated by the 
Davey-Stewartson hierarchy of $2+1$-dimensional soliton equations. These 
deformations of surfaces inherit all remarkable properties of the soliton 
equations. Geometrically, such deformations are characterized by the 
invariance of an infinite set of functionals over surfaces. The simplest of 
them is given by the Willmore functional.

Note that the Weierstrass type representations for particular classes of 
surfaces in $R^4$ have been discussed in \bib{r25},\bib{r28} 
and recently in 
\bib{r29}. Formulae of the Weierstrass type for immersion of surfaces into 
$R^4$ can be derived also within the quaternionic approach \bib{r30}.

  \section{Generalized Weierstrass representations}
\setcounter{equation}{0}

Main steps in our construction of the generalized Weierstrass formulae for 
surfaces in four-dimensional spaces are basically the same as those in 
$\RR^3$ (see also \bib{r29}).

\begin{proposition}

The generalized Weierstrass formulae 
\begin{eqnarray}
X^1 + i X^2 & = &  \int_{\Gamma} { \left( -
\varphi_1 \varphi_2 dz' +
\psi_1 \psi_2  d\zba' \right) }
\; , \nonumber \\
X^1 - i X^2 & = &  \int_{\Gamma} { \left(
\psiba_1 \psiba_2 dz' - \vphiba_1 \vphiba_2  d\zba' \right)}
\; , \nonumber \\
X^3 + i X^4 & = &  \int_{\Gamma} { \left(
\psiba_2 \varphi_1 dz'+ \psi_1 \vphiba_2 d\zba' \right)}
\; , \nonumber \\
X^3 - i X^4 & = & \int_{\Gamma} { \left(
\psiba_1 \varphi_2 dz' + \psi_2 \vphiba_1 d\zba' \right) }
\label{eq:2.1}
\end{eqnarray}
where
\begin{eqnarray}
\begin{array}{l}
\psi_{1 z}=p \varphi_1 \;\; ,\\
\varphi_{1 \zba}=-\pba \psi_1
\end{array}
&   & 
\begin{array}{l}
\psi_{2 z}=\pba \varphi_2 \; , \\
\varphi_{2 \zba}=-p \psi_2 \;\;\; ,
\end{array}
\label{eq:2.2}
\end{eqnarray}
$\psi_{\alpha}$, $\varphi_{\alpha}$ ($\alpha=1,2$), $p$ are complex-valued
functions of $z$, $\zba$, $\Gamma$ is a contour in $\CC$, define the 
conformal immersion of a surface into $\RR^4$: $X^i(z,\zba): 
\CC \rightarrow \RR^4$. The induced metric of a surface is of the form 
\beq
ds^2=u_1 u_2 dz d\zba
\label{eq:2.3}
\ee
where
$u_{\alpha}=|\psi_{\alpha}|^2+|\varphi_{\alpha}|^2 $, $\alpha=1,2$, the Gaussian and 
mean curvatures are respectively 
\begin{eqnarray}
K=-\frac{2}{u_1 u_2}
\left[ \log { \left( u_1 u_2 \right) } \right] _{z \zba}
\;\; , \;\; \vec{H}^2=\frac{4|p|^2}{u_1 u_2}
\label{eq:2.4}
\end{eqnarray}
the total squared mean curvature 
$W=\int {\vec{H}^2  \left[ dS \right] }$
(Willmore functional) is
\beq
W=4 \int{ |p|^2 dx \;dy} \;\;.
\label{eq:2.5}
\ee

\end{proposition}

Indeed, equations (\ref{eq:2.2}) imply that 
\begin{eqnarray}
\left( \psi_1 \psi_2 \right)_{z} = - \left( \varphi_1 \varphi_2
\right)_{\zba}
\;\; , \;\;
\left( \psi_1 \vphiba_2 \right)_{z}=\left( \varphi_1 \psiba_2
\right)_{\zba} \;\;.
\label{eq:2.6}
\end{eqnarray}
In virtue of (\ref{eq:2.6}) the {\em r.h.s.} of (\ref{eq:2.1})
do not depend on the contour $\Gamma$ and, hence, the coordinates $X^i$ ($i
=1,2,3,4$) given by (\ref{eq:2.1}) are defined uniquely up to the 
displacement constants. Induced metric and Gaussian curvature are 
calculated straightforwardly. 
The mean curvature vector $\vec{H}=\frac{\vec{X}_{z\zba}}{g_{z\zba}}$ 
is given by
\begin{eqnarray}
\vec{H} & =  & \frac{2}{u_1 u_2} 
\left[ Re\left( p \psi_2 \varphi_1 + \pba \psi_1 \varphi_2 \right), 
Im\left(p \psi_2 \varphi_1 +\pba \psi_1 \varphi_2 \right) , \right. 
\nonumber \\
& & \left. \;\;\;\;\;\;\;\;\;\;\; \;\;
Re\left( p \varphi_1 \vphiba_2 - 
\pba \psi_1 \psiba_2 \right), 
Im\left( p \psiba_1 \psi_2 -\pba \vphiba_1 \varphi_2 \right)
 \right] \; 
\label{eq:2.7}
\end{eqnarray}
where $Re$, $Im$ denote the real and imaginary parts, respectively.
The formula (\ref{eq:2.7}) yields the expression (\ref{eq:2.4}) for 
$\vec{H}^2$ and, finally, (\ref{eq:2.5}).

Since equations (2.2) contain two arbitrary real-valued functions of two 
variables ($Re(p)$ and $Im(p)$) the Weierstrass representation
(\ref{eq:2.1})-(\ref{eq:2.2}) provides us a generic surface in $\RR^4$.

\begin{corollary}
 Surfaces of a constant mean curvature $H$ ($\vec{H}^2=H^2=constant$) in 
$\RR^4$ are generated by the formulae (\ref{eq:2.1}) where $\psi_{\alpha}$, 
$\varphi_{\alpha}$ ($\alpha=1,2$) obey the system of equations 
\begin{eqnarray}
\psi_{\alpha z} & = & \frac{H}{2} e^{i\theta}
\sqrt{\left(|\psi_1|^2+|\varphi_1|^2 \right)
\left(|\psi_2|^2+|\varphi_2|^2 \right)} \; \varphi_\alpha
\;\; , \nonumber \\
\varphi_{\alpha \zba} & = & -\frac{H}{2} e^{-i\theta}
\sqrt{\left(|\psi_1|^2+|\varphi_1|^2 \right)
\left(|\psi_2|^2+|\varphi_2|^2 \right)} \; \psi_\alpha
\;\; , \;\; \alpha=1,2 
\label{eq:2.8}
\end{eqnarray}
where $\theta(z,\zba)$ is an arbitrary function. 
\end{corollary}

 At the particular case $\pba=p$ the formulae 
(\ref{eq:2.1})-(\ref{eq:2.5}) are reduced to those derived in \bib{r29} 
with the substitution $X^1 \leftrightarrow X^2$, 
$X^3 \leftrightarrow -X^3$. This case corresponds to a special class of 
surfaces in $\RR^4$. The reduction $\psi_1=\pm \psi_2$, $\varphi_1=
\pm \varphi_2$ (consequently $\pba=p$) converts 
(\ref{eq:2.1})-(\ref{eq:2.5}) into the generalized Weierstrass 
representation (\ref{eq:1.1})-(\ref{eq:1.2}) for generic surfaces in
$\RR^3$. 

Note that a linear system of the form (\ref{eq:2.2}) arises also as the
restriction of the Dirac equation to a surface in $\RR^4$ \bib{r31}.

The case of the pseudo-Euclidean space $\RR^{2,2}$ with the metric
$g_{ik}=diag(1,$ $1,-1,-1)$ is rather similar to that of $\RR^4$. 

\begin{proposition}
The generalized Weierstrass formulae
\begin{eqnarray}
X^1 + i X^2 & = &  \int_{\Gamma} { \left( 
\varphi_1 \varphi_2 dz' +
\psi_1 \psi_2 d\zba' \right) }
\; , \nonumber \\
X^1 - i X^2 & = &   \int_{\Gamma} { \left(
\psiba_1 \psiba_2 dz' + \vphiba_1 \vphiba_2  d\zba' \right)}
\; , \nonumber \\
X^3 + i X^4 & = & i \int_{\Gamma} { \left(
\psiba_1 \varphi_2 dz'+ \psi_2 \vphiba_1 d\zba' \right)}
\; , \nonumber \\
X^3 - i X^4 & = & -i\int_{\Gamma} {\left(
\psiba_2 \varphi_1 dz' + \psi_1 \vphiba_2 d\zba' \right) }
\label{eq:2.9}
\end{eqnarray}
where
\begin{eqnarray}
\begin{array}{l}
\psi_{1 z}=p \varphi_1 \;\; ,\\
\varphi_{1 \zba}=\pba \psi_1
\end{array}
&   &
\begin{array}{l}
\psi_{2 z}=\pba \varphi_2 \; , \\
\varphi_{2 \zba}=p \psi_2 \;\;\; ,
\end{array}
\label{eq:2.10}
\end{eqnarray}
$\psi_{\alpha}$, $\varphi_{\alpha}$, $p$ are complex-valued functions,
$\Gamma$ is a contour in $\CC$, define the conformal immersion
$\vec{X}:\CC \rightarrow \RR^{2,2}$ of a surface into the space 
$\RR^{2,2}$. The induced metric is 
\beq
ds^2=v_1 v_2 dz d\zba
\label{eq:2.11}
\ee
where $v_{\alpha}=|\psi_{\alpha}|^2-|\varphi_{\alpha}|^2$, 
$\alpha=1,2$, the Gaussian and mean curvature are of the form
\begin{eqnarray}
K=-\frac{2}{v_1 v_2}
\left[ \log { \left( v_1 v_2 \right) } \right] _{z \zba}
\;\; , \;\; \vec{H}^2=-\frac{4|p|^2}{v_1 v_2} \;\;
\label{eq:2.12}
\end{eqnarray}
and the Willmore functional $W=\int{\vec{H}^2 \left[ dS
\right] }$ is given by
\beq
W=-4\int {|p|^2 dx dy} \;\;.
\label{eq:2.13}
\ee

\end{proposition}

The proof is similar to the case of $\RR^4$, only now the equations
(\ref{eq:2.10}) give $\left( \psi_1 \psi_2 \right)_z=
\left( \varphi_1 \varphi_2 \right)_{\zba}$ and 
$\left( \psi_1 \vphiba_2 \right)_z=
\left( \varphi_1 \psiba_2 \right)_{\zba}$. In the particular case
$\pba=p$ the formulae (\ref{eq:2.9})-(\ref{eq:2.13}) are reduced to those 
obtained in \bib{r29} with the substitution
$X^1 \leftrightarrow X^2$, $X^3 \leftrightarrow X^4$. In contrast to
\bib{r29} the formulae (\ref{eq:2.9})-(\ref{eq:2.10}) allow to represent
an arbitrary surface in $\RR^{2,2}$.

Surfaces in $\RR^{2,2}$ with the constant $\vec{H}^2$ are generated by
(\ref{eq:2.9}) where $\psi_{\alpha}$, $\varphi_{\alpha}$ obey the system 
(\ref{eq:2.8}) with obvious change of signs. Note that the space
$\RR^{2,2}$ is of importance also in string theory \bib{r32}.

Conformal immersions into the Minkowski space $\RR^{3,1}$ are given by 
slightly different formulae.

\begin{proposition}
The Weierstrass type formulae
\begin{eqnarray}
X^1+iX^2 & = & \int_{\Gamma} { \left(
\psiba_2 \varphi_1  dz'+
\psi_1 \vphiba_2 d\zba' \right)  }
\; , \nonumber \\
X^1-iX^2 & = & \int_{\Gamma} { \left(
\psiba_1 \varphi_2  dz'+ \psi_2 \vphiba_1 d\zba' \right) }
\; , \nonumber \\
X^3+X^4 & = & \int_{\Gamma} { \left(
\psiba_1 \varphi_1  dz'+ \psi_1 \vphiba_1 d\zba' \right) }
\; , \nonumber \\
X^3-X^4 & = & -\int_{\Gamma} { \left(
\psiba_2 \varphi_2  dz'+ \psi_2 \vphiba_2 d\zba' \right) }
\label{eq:2.14}
\end{eqnarray}
where 
\begin{eqnarray}
\begin{array}{l}
\psi_{\alpha z}=p \varphi_{\alpha} \; , \\
\varphi_{\alpha \zba}=q \psi_\alpha
\end{array}
 &, & \alpha=1,2
\label{eq:2.15}
\end{eqnarray}
$p$ and $q$ are real-valued functions, $\Gamma$ is a contour in $\CC$, 
define the conformal immersion of a surface into the Minkowski space
$\vec{X}:\CC \rightarrow \RR^{3,1}$. The induced metric on a surface is
\beq
ds^2=|\psi_1 \varphi_2-\psi_2 \varphi_1|^2 dz d\zba \;\; ,
\label{eq:2.16}
\ee
the mean curvature $\vec{H}^2$ and the Willmore functional are given 
respectively by
\beq
\vec{H}^2=-\frac{4pq}{|\psi_1\varphi_2-\psi_2\varphi_1|^2}
\;\; , \;\; W=-4 \int {pq dx dy} .
\label{eq:2.17}
\ee

\end{proposition}

In this case the linear system (\ref{eq:2.15}) implies that
$$ \left( \psi_{\alpha} \vphiba_{\beta} \right)_z=
\left( \varphi_{\alpha} \psiba_{\beta} \right)_{\zba} \;\;\;\;
\alpha,\beta=1,2 $$
that guarantee an independence of the {\em r.h.s.} of (\ref{eq:2.14}) on the 
choice of the contour $\Gamma$ of integration. The rest is straightforward.

Since again one has two arbitrary real-valued functions $p$ and $q$, the 
Weierstrass type formulae (\ref{eq:2.14}), (\ref{eq:2.15}) allow us to 
construct any surface immersed into $\RR^{3,1}$.

Surfaces of constant $\vec{H}^2$ are generated by the formulae 
(\ref{eq:2.14}), where 
$\psi_{\alpha}$ and $\varphi_{\alpha}$
obey the system of equations (\ref{eq:2.15}) with 
the constraint 
\beq
pq=-\frac{1}{4} \vec{H}^2 |\psi_1 \varphi_2 -\psi_2 \varphi_1|^2 \;\;.
\label{eq:2.18}
\ee
Special classes of surfaces which correspond to the cases
$q=constant$ or $p=\pm q$ could be of particular interest.

Differential versions of all three generalized Weierstrass representations 
given above can be written in the following common form
\beq
d \left( \sum_{i=1}^4 { \tau_i X^i} \right) =\Phi^{\dagger}_{2} 
\left(
\begin{array}{cc}
0 & dz \\
d\zba & 0
\end{array}
\right)
\Phi_1
\label{eq:2.19}
\ee
where $\dagger$ denotes Hermitian conjugation. In the case of 
immersion into $\RR^4$ one has
$$
\tau_1=\sigma_1 \;\;,\;\;
\tau_2=\sigma_2 \;\;,\;\;
\tau_3=\sigma_3 \;\;,\;\;
\tau_4=i\sigma_4
$$
and
\beq
\Phi_{\alpha}= \left(
\begin{array}{cc}
\psi_{\alpha} & -\vphiba_{\alpha} \\
\varphi_{\alpha} & \psiba_{\alpha}
\end{array}
\right) \; \; , \;\; \alpha=1,2
\label{eq:2.20}
\ee
where $\sigma_1$, $\sigma_2$, $\sigma_3$ are the standard Pauli matrices 
and $\sigma_4$ is an identical $2 \times 2$ matrix. At the $\RR^{2,2}$ case
$$
\tau_1=\sigma_1 \;\;,\;\;
\tau_2=\sigma_2 \;\;,\;\;
\tau_3=i\sigma_3 \;\;,\;\;
\tau_4=\sigma_4
$$
and
\beq
\Phi_{\alpha}= \left(
\begin{array}{cc}
\psi_{\alpha} & \vphiba_{\alpha} \\
\varphi_{\alpha} & \psiba_{\alpha}
\end{array}
\right) \; \; , \;\; \alpha=1,2
\label{eq:2.21}
\ee
Finally, the immersion into the Minkowski space $\RR^{3,1}$ corresponds to 
$$
\tau_i=\sigma_i \;\; \;\; (i=1,2,3,4)
$$
and
\beq
\Phi_1=\Phi_{2}= \left(
\begin{array}{cc}
\psi_{1} & \psi_{2} \\
\varphi_{1} & \varphi_{2}
\end{array}
\right) 
\; \; .
\label{eq:2.22}
\ee
In fact, one can start with the formulae (\ref{eq:2.19}) to derive the 
Weierstrass representations in the forms (\ref{eq:2.1})-(\ref{eq:2.2}),
(\ref{eq:2.9})-(\ref{eq:2.10}) and (\ref{eq:2.14})-(\ref{eq:2.15}). Indeed, 
one can show that the $1-$form in the {\em r.h.s.} of (\ref{eq:2.19})
is closed if the $2 \times 2$ matrices $\Phi_1$, $\Phi_2$ obey the 
Dirac equations
\beq
\left(
\begin{array}{cc}
\partial_{z} & 0 \\
0 & \partial_{\zba}
\end{array}
\right) \Phi_1 =
\left(
\begin{array}{cc}
u & p \\
q & v
\end{array}
\right) \Phi_1 \;\; , \;\;
\left(
\begin{array}{cc}
\partial_{z} & 0 \\
0 & \partial_{\zba} 
\end{array}
\right) \Phi_2= 
\left(
\begin{array}{cc}
-\overline{v} & \overline{p} \\
\overline{q} & -\overline{u}
\end{array}
\right) \Phi_2
\label{eq:2.23}
\ee
where $p$, $q$, $u$, $v$ are arbitrary complex-valued functions. Functions 
$u$ and $v$ always can be converted to zeros by
gauge transformation (redefinition of $\Phi$).
Then the reality conditions for $X^i$ are 
satisfied if matrices $\Phi_{\alpha}$ have the form (\ref{eq:2.20}), 
(\ref{eq:2.21}) or (\ref{eq:2.22}) while the functions $p$, $q$ should obey 
the constraints $p+\overline{q}=0$, $p-\overline{q}=0$ and
$p=\overline{p}$, $q=\overline{q}$, respectively. Consequently, the 
corresponding formula (\ref{eq:2.19}) gives rise to the Weierstrass 
representations considered above. 

A formula of the type (\ref{eq:2.19}) appears naturally \bib{r30} in the 
quaternionic approach to surfaces (see also \bib{r33}-\bib{r35})
which could 
provide an invariant formulation of the construction presented above.

Following to \bib{r29}, one can extend the generalized Weierstrass 
representations to generic surfaces immersed into four-dimensional Riemann 
spaces. In particular, in the cases of spaces of constant curvature and 
conformally-flat spaces the Willmore functional still has the form
(\ref{eq:2.5}), (\ref{eq:2.13}) and (\ref{eq:2.17}). 

Explicit construction of surfaces via the Weierstrass formulae requires the 
resolving the linear systems (\ref{eq:2.2}), (\ref{eq:2.10})
and (\ref{eq:2.15}). The methods of solving these type of problems are 
well-developed now within the theory of $(2+1)-$dimensional integrable
(soliton) equations (see {\em e.g} \bib{r36}-\bib{r38}). Using these methods one 
can construct broad classes of surfaces explicitly. The results will be 
presented in a separate paper.

\section{Integrable deformations}
\setcounter{equation}{0}

Now, following the general approach of \bib{r10}-\bib{r11}, we will 
construct integrable deformations of surfaces generated by the Weierstrass 
type formulae. The basic idea is to use deformations of the functions $p$, 
$q$, $\psi$, $\varphi$ compatible with the linear equations (\ref{eq:2.2}), 
(\ref{eq:2.10}) and (\ref{eq:2.15}).

All of them are the particular cases of the linear system 
\beq
\begin{array}{l}
\psi_{z}=p \varphi \;\; ,\\
\varphi_{\zba}=q \psi
\end{array}
\label{eq:3.1}
\ee
where $p$ and $q$ are in general complex-valued functions. In soliton 
theory this system is known as the Davey-Stewartson II (DSII) linear 
problem (see {\em e.g.} \bib{r36}-\bib{r38}).

An infinite hierarchy of nonlinear differential equations associated with
(\ref{eq:3.1}) is referred as the DSII hierarchy. It arises as the 
compatibility conditions of (\ref{eq:3.1}) with the systems 
\bib{r36}-\bib{r38}
\beq
\begin{array}{l}
\psi_{t_n}=A_n \psi+B_n \varphi \; \; , \\
\varphi_{t_n}=C_n \psi+D_n \varphi \;\; .
\end{array}
\label{eq:3.2}
\ee
where $t_n$
are new (deformation) variables and $A_n$, $B_n$, $C_n$,
$D_n$ are differential operators of $n-$th order. At $n=1$ one gets
the linear system
\beq
\begin{array}{l}
p_{t_1}=\alpha p_{\zba} +\gamma p_z\; \; , \\
q_{t_1}=\gamma q_{z} +\alpha q_{\zba}
\end{array}
\label{eq:3.3}
\ee
where $\alpha$, $\gamma$ are arbitrary constants. The corresponding 
operators in (\ref{eq:3.2}) are
\beq
A_1=\alpha \partial_{\zba} \;,\;B_1=\gamma p\;,\; C_1=\alpha q\;,\; 
D_1=\gamma \partial_z
 \;\;\;.
\label{eq:3.4}
\ee
Higher equations are nonlinear ones. At $n=2$ one has the system
\bib{r36}-\bib{r38}
\beq
\begin{array}{l}
p_{t_2}=\alpha_2 \left( p_{zz} + p_{\zba \zba} +up \right) \; \; , \\
q_{t_2}=-\alpha_2 \left( q_{zz} + q_{\zba \zba} +uq \right) \; \; , \\
u_{z \zba}=-2(pq)_{zz}-2(pq)_{\zba \zba}
\end{array}
\label{eq:3.5}
\ee 
where $\alpha_2$ is an arbitrary constant. For the system (\ref{eq:3.5})
\begin{eqnarray}
\begin{array}{l}
A_{2}=\alpha_2 \left( \partial_{\zba}^{2}+w_1 \right)
\; ,\\
C_{2}=-\alpha_2 \left( q_{\zba} - q \partial_{\zba} \right) \; ,
\end{array}
&  &
\begin{array}{l}
B_{2}=\alpha_2 \left( p_z- p\partial_z \right)
\; ,\\
D_{2}=-\alpha_2 \left( \partial_{z}^{2}+w_2 \right)
\end{array}
\label{eq:3.6}
\end{eqnarray}
where 
$$ w_{1z}=-2(pq)_{\zba}\;\;,\;\; w_{2 \zba}=-2(pq)_{z}\;\;,\;\;
u=w_1-w_2\;\;.$$
For the $t_3-$deformations one has (see {\em e.g.} \bib{r37})
\beq
\begin{array}{lll}
p_{t_3} & = & p_{zzz}+p_{\zba \zba \zba} +
3p_z \partial^{-1}_{\zba} (pq)_{z} 
+ 3p_{\zba} \partial^{-1}_{z}(pq)_{\zba}
+3p \partial^{-1}_{\zba}(qp_z)_{z} 
+ 3p \partial^{-1}_{z} (qp_{\zba})_{\zba} 
\; , \nonumber \\
q_{t_3} & = & q_{zzz}+q_{\zba \zba \zba} +
3q_z \partial^{-1}_{\zba} (pq)_{z}
+ 3q_{\zba} \partial^{-1}_{z}(pq)_{\zba}
+3q \partial^{-1}_{\zba}(pq_z)_{z}
+ 3q \partial^{-1}_{z} (pq_{\zba})_{\zba} \;\;.
\end{array}
\label{eq:3.7}
\ee
In this case
\begin{eqnarray}
\begin{array}{l}
A_3=\partial_{\zba}^3+3 \left[ \partial_{z}^{-1} (pq)_{\zba}
\right] \partial_{\zba} +3 \left[ \partial_z^{-1} 
(qp_{\zba})_{\zba} \right] \; ,\\
B_3=-p \partial^2_{z} + p_z \partial_z-p_{zz}-
3p \left[ \partial_{\zba}^{-1}(pq)_{z} \right] \; ,\\
C_3=-q \partial^2_{\zba} + q_{\zba} \partial_{\zba}-q_{\zba \zba}
-3q \left[ \partial^{-1}_z (pq)_{\zba} \right] \; ,\\
D_3=\partial_{z}^3+3 \left[ \partial_{\zba}^{-1} (pq)_{z}
\right] \partial_{z} +3 \left[ \partial_{\zba}^{-1}
(pq_{z})_{z} \right] \; .
\end{array}
\label{eq:3.8}
\end{eqnarray}
All equations of this DSII hierarchy, known also as the
$2-$component KP 
hierarchy, are the $(2+1)-$dimensional soliton equations. They are 
integrable by the inverse spectral transform method
(see {\em e.g.} \bib{r36}-\bib{r38}) and have all remarkable properties typical 
for soliton equations: namely, they possess infinite classes of explicit 
solutions, including multi-soliton solutions, infinite symmetry algebra, 
Darboux and Backlund transformations etc. They have an infinite set of 
integrals of motion. The simplest of them is
\beq
C_1=\int {p q \; dx dy} 
\label{eq:3.9}
\ee
while other integrals of motion are non-local. Note that (\ref{eq:3.9})
is the integral of motion for the whole DSII hierarchy ($C_{1t_n}=0$). 

So, now we assume that all quantities (except $z$, $\zba$) in
the systems
(\ref{eq:2.2}), (\ref{eq:2.10}) and (\ref{eq:2.15}) and, consequently, 
coordinates $X^i$ depend on the deformation parameters $t_n$ ($n=1,2,
\ldots$) (times). Then, we assume that this dependence on $t_n$ is such 
that there are operators $A_n$, $B_n$, $C_n$, $D_n$ such 
that (\ref{eq:3.2}) holds. The compatibility conditions of
(\ref{eq:2.2}), (\ref{eq:2.10}) or (\ref{eq:2.15}) with
(\ref{eq:3.2}) fix the dependence of $\psi$, $\varphi$ and $p$, $q$ on 
$t_n$ and, 
consequently, define the deformations of surfaces. Concrete cases are 
governed by different specialization (reductions) of the DSII 
hierarchy.

Let us consider first immersions into the Minkowski space $\RR^{3,1}$. In 
this case $p$ and $q$ are real-valued functions. The corresponding 
deformations are generated by the "real" DSII hierarchy (with real-
valued $p$ and $q$). In particular, in equations (\ref{eq:3.3}),
(\ref{eq:3.5}) the constants $\alpha$, $\gamma$, $\alpha_2$ are real.
Since (\ref{eq:3.9}) is obviously the integral of motion also for
real-valued $p$ and $q$, then, in virtue of (\ref{eq:2.17}), the Willmore 
functional $W$ remains invariant under these DSII deformations. 

For the Weierstrass representations in the spaces $\RR^4$ and $\RR^{2,2}$ 
we have the linear system (\ref{eq:3.1}) with the following reductions:
\beq
\left(
\begin{array}{cc}
0 & p_1 \\
q_1 & 0
\end{array}
\right)=
\left(
\begin{array}{cc}
0 & p \\
\varepsilon \pba & 0
\end{array}
\right)
\;\;\;\; for \;\;\; \Phi_1
\label{eq:3.10}
\ee
and
\beq
\left(
\begin{array}{cc}
0 & p_2 \\
q_2 & 0
\end{array}
\right)
=\left(
\begin{array}{cc}
0 & \pba \\
\varepsilon p & 0
\end{array}
\right)
\;\;\;\; for \;\;\; \Phi_2
\label{eq:3.11}
\ee
where $\varepsilon=-1$ for $\RR^4$ and $\varepsilon=1$ for $\RR^{2,2}$. 
Both the reductions (\ref{eq:3.10}) and (\ref{eq:3.11}) are admissible by 
all equations of the DSII hierarchy if one chooses (see {\em e.g.} 
\bib{r36}-\bib{r38})
\begin{eqnarray}
\begin{array}{l}
A_{2n-1}=\partial_{\zba}^{2n-1}+\ldots, \\ 
A_{2n}=\pm i\partial_{\zba}^{2n}+\ldots,
\end{array}
&  &
\begin{array}{l}
D_{2n-1}=\partial_{z}^{2n-1}+\ldots, \\
D_{2n}=\mp i\partial_{z}^{2n}+\ldots \;\;\;\;\;\;(n=1,2,3,\dots).
\end{array}
\nonumber
\end{eqnarray}
In our case we have different linear problems for $\psi_1$, $\varphi_1$
and $\psi_2$, $\varphi_2$. To have the same equations for $p$ it is 
sufficient to take
\begin{eqnarray}
\begin{array}{l}
A_{2n-1}=\partial_{\zba}^{2n-1}+\ldots, \\
A_{2n}=i\partial_{\zba}^{2n}+\ldots,
\end{array}
&  &
\begin{array}{l}
D_{2n-1}=\partial_{z}^{2n-1}+\ldots, \\
D_{2n}=-i\partial_{z}^{2n}+\ldots
\end{array}
\label{eq:3.12}
\end{eqnarray}
in equations (\ref{eq:3.2}) for $\psi_1$, $\varphi_1$ and
\begin{eqnarray}
\begin{array}{l}
A_{2n-1}=\partial_{\zba}^{2n-1}+\ldots, \\
A_{2n}=-i\partial_{\zba}^{2n}+\ldots,
\end{array}
&  &
\begin{array}{l}
D_{2n-1}=\partial_{z}^{2n-1}+\ldots, \\
D_{2n}=i\partial_{z}^{2n}+\ldots
\end{array}
\label{eq:3.13}
\end{eqnarray}
in the case of $\psi_2$, $\varphi_2$. In particular, one gets
\beq
\begin{array}{l}
p_{t_2}=i \left( p_{zz}+p_{\zba \zba} +u p \right)
\; \; , \\
u_{z \zba}=-2 \varepsilon |p|^2_{zz}-2\varepsilon |p|^2_{\zba \zba}
\end{array}
\label{eq:3.14}
\ee
and
\beq
\begin{array}{l}
\psi_{1t_2}=i \left( \partial_{\zba}^2+w_1 \right) \psi_1 +
i \left( p_{z}-p \partial_{z} \right) \varphi_1
\; \; , \\
\varphi_{1t_2}=-i \varepsilon \left( \pba_{\zba}-\pba \partial_{\zba} \right)
\psi_1 -i \left( \partial_{z}^2+w_2 \right) \varphi_1
\end{array}
\label{eq:3.15}
\ee
while
\beq
\begin{array}{l}
\psi_{2t_2}=-i \left( \partial_{\zba}^2+w_1 \right) \psi_2 -
i \left( \pba_{z}-\pba \partial_{z} \right) \varphi_2
\; \; , \\
\varphi_{2t_2}=i \varepsilon \left( p_{\zba}-p\partial_{\zba} \right)
\psi_2+i\left( \partial_{z}^2+w_2 \right) \varphi_2
\end{array}
\label{eq:3.16}
\ee
where $w_{1z}=-2\varepsilon |p|^2_{\zba}$, 
$w_{1\zba}=2\varepsilon |p|^2_{z}$. 

The $t_3$ deformation is given now by the equation
\beq
\begin{array}{lll}
p_t& = & p_{zzz}+p_{\zba \zba \zba} +3 \varepsilon p_z \partial^{-1}_{\zba}
(|p|^2_z) +3 \varepsilon p_{\zba} \partial^{-1}_{z} (|p|^2_{\zba})+ \\
 & & +3 \varepsilon p \partial^{-1}_{\zba}(\pba p_z)_z+
3 \varepsilon p \partial^{-1}_{z}(\pba p_{\zba})_{\zba}
\end{array}
\label{eq:3.17}
\ee
and the deformations of  $\psi_1$, $\varphi_1$ and
 $\psi_2$, $\varphi_2$ are given by (\ref{eq:3.2}), (\ref{eq:3.8}) with the 
reduction (\ref{eq:3.10}) and (\ref{eq:3.11}), respectively.

Thus, in the cases of $\RR^4$ and $\RR^{2,2}$ deformations of surfaces are 
generated by the proper DSII equation (\ref{eq:3.14}) and the 
corresponding hierarchy. Properties of solutions of the DSII equation
(\ref{eq:3.14}) are essentially different for different signs of 
$\varepsilon$. Consequently, the properties of deformations of surfaces in
$\RR^4$ and $\RR^{2,2}$ will differ too. 

In both cases $C_1=\int {|p|^2 dxdy}$ is the integral of motion for 
the whole hierarchy. Hence, the Willmore functional $W$ for surfaces 
immersed into $\RR^4$ and $\RR^{2,2}$ is invariant under deformations 
generated by the DSII hierarchy.

Thus, though the deformations for surfaces in $\RR^4$, $\RR^{3,1}$ and
$\RR^{2,2}$ are governed by different nonlinear integrable equations, they 
have the following common property.

\begin{proposition}
The DSII hierarchy generates integrable deformations of surfaces immersed 
into $\RR^4$, $\RR^{3,1}$ and $\RR^{2,2}$ via the generalized Weierstrass
representations. The Willmore functionals W for surfaces in
$\RR^4$, $\RR^{3,1}$ and $\RR^{2,2}$ are invariant under the
corresponding deformations ($W_{t_n}=0$).
\end{proposition}

DSII hierarchy of integrable equations is well studied (see
{\em e.g.} \bib{r36}-\bib{r38}). This provides us a broad class of deformations 
of surfaces in $\RR^4$, $\RR^{3,1}$ and $\RR^{2,2}$ given explicitly. 
Moreover, since the inverse spectral transform method allows us to 
linearize the initial-value problem 
$${p(z,\zba,t_n=0), q(z,\zba,t_n=0)}
\rightarrow {p(z,\zba,t_n), q(z,\zba,t_n)}$$
for soliton equations of the 
DSII hierarchy (see {\em e.g.} \bib{r36}-\bib{r38}), then the generalized 
Weierstrass formulae allows us to linearize the initial-value problem for 
the deformation of surfaces $\vec{X}(z,\zba,0) \rightarrow 
\vec{X}(z,\zba,t_n)$. In virtue of all that, the deformations generated by 
the DSII hierarchy can be referred as integrable one.

Higher integrals of motion for the DSII hierarchy are also certain 
functionals on surfaces invariant under deformations generated by the
DSII hierarchy.
Since the Willmore functional $W$ is invariant under the 
conformal transformations in four-dimensional spaces, then it is quite 
natural to suggest

\begin{conjecture}
Higher integrals of motion for the DSII hierarchy are functionals on 
surfaces in $\RR^4$, $\RR^{3,1}$ and $\RR^{2,2}$ which are invariant under 
conformal transformations in these spaces.
\end{conjecture}

For surfaces in $\RR^3$ an analogous conjecture has been proved in 
\bib{r17}.

In the particular case $p=\pba$ the Weierstrass representations for $\RR^4$ 
and $\RR^{2,2}$ are reduced to those derived in \bib{r29}. The reduction 
$p=\pba$ is admissible only by deformations associated with odd times 
$t_{2n-1}$. In this case equation (\ref{eq:3.17}) is nothing but the 
modified Veselov-Novikov equation and the DSII sub-hierarchy associated 
with only odd times is converted into the mVN hierarchy. This specialized 
class of surfaces in $\RR^4$ and $\RR^{2,2}$ and their deformations via the 
mVN hierarchy have been considered in \bib{r29}.

For the case of Minkowski space $\RR^{3,1}$ there is a class of surfaces 
for which $q=1$. Integrable deformations of such class of surfaces are 
generated by the DSII hierarchy under the reduction $q=1$.
This reduction is compatible only with odd times 
deformations. For instance, the $t_3-$deformation is of the form
(reduction $q=1$ of the system (\ref{eq:3.7}))
\beq
p_{t}= p_{zzz}+ p_{\zba \zba \zba}+3 \left[ p \partial_{\zba}^{-1}
(p_z) \right]_z + 3 \left[ p \partial_{z}^{-1}
(p_{\zba}) \right]_{\zba}
\label{eq:3.18}
\ee
that is the Veselov-Novikov equation \bib{r36}-\bib{r38}. So, integrable 
deformations of surfaces with $q=1$ are generated by the Veselov-Novikov 
hierarchy.

Generalized Weierstrass formulae with specialized $p=p(x)$ and $q=q(x)$ 
($x=Re(z)$) give rise obviously to surfaces of "revolution" in four 
dimensional spaces. In this case the DSII hierarchy is reduced to the 
AKNS hierarchy of the $(1+1)-$dimensional soliton equations
(see {\em e.g.} \bib{r36}-\bib{r38}). Equation
(\ref{eq:3.14}) is converted into the well-known nonlinear Schroedinger
(NLS) equation
\beq
ip_{t_2}+\frac{1}{2} p_{xx}-2\varepsilon |p|^2p=0 \;\;.
\label{eq:3.19}
\ee
So, integrable deformations of surfaces of revolution in $\RR^4$ and
$\RR^{2,2}$ are generated by the NLS hierarchy. The NLS equation
arises also in the description of integrable motion of space
curves in $\RR^3$ \bib{r39}. So, there is an intimate
connection between integrable deformations of 
surfaces of revolution in $\RR^4$ and integrable motions of space curves in 
$\RR^3$.

Finally, we note that in the particular case $\psi_1=\pm \psi_2$, 
$\varphi_1=\pm \varphi_2$ the Weierstrass representation 
(\ref{eq:2.1})-(\ref{eq:2.2}) is reduced to that in $\RR^3$ 
\bib{r10}-\bib{r11}. Since then $\pba=p$, the DSII
hierarchy is reduced to 
the mVN hierarchy which generates integrable deformations of surfaces in 
$\RR^3$ \bib{r10}-\bib{r11}.

\hfill
\hfill

{\bf Acknowledgements.} The author is grateful to F. Pedit and U. Pinkall 
for very useful and stimulating discussions. He thanks the 
Sonderforschungsbereich 288 for the kind hospitality and financial support 
during the visit to the Technical University of Berlin where the main part 
of this work was done.


\newpage

   \begin{centerline} 
   {\bf REFERENCES}
   \end{centerline}

   \begin{enumerate}

  \item \label{r1}
	Darboux G., {\it Lecons sur la th\'{e}orie des surfaces et 
	les applications geometriques du calcul infinitesimal}, 
	{\em t. 1-4}, Gauthier-Villars, Paris, 1877-1896.	
  \item \label{r2}
	Bianchi L., {\it Lezioni di Geometria Differenziale}, 
	{\it 2nd ed.}, Spoerri, Pisa, 1902.
	
  \item \label{r3}
	Eisenhart L.P., {\it A treatise on the differential geometry
	of Curves and Surfaces}, Dover, New York, 1909.

  \item \label{r4}
	Yau S.T. (Ed.), {\it Seminar on Differential Geometry}, 
	Princeton Univ., Princeton, 1992.
  \item \label{r5}
	Chern S.S., in: {\it Differential Geometry and Complex Analysis} 
	(Chavel I., Farkas H.M., Eds.), Springer, Berlin, 1985.
  \item \label{r6}
	Nelson D., Piran T. and Weinberg S. (Eds.),
	{\it Statistical mechanics 
	of membranes and surfaces}, Worls Scientific, Singapore, 1989.
 
  \item \label{r7}
	Gross D.J., Piran T. and Weinberg S. (Eds.), {\it Two 
	dimensional quantum gravity and random surfaces}, 
	World Scientific, Singapore, 1992. 
	
  \item \label{r8}
	David F., Ginsparg P. and Zinn-Justin Y. (Eds.),
	{\it Fluctuating geometries in Statistical Mechanics 
	and Field Theory}, Elsevier Science, Amsterdam, 1996.
  \item \label{r9} 
	Kenmotsu K., {\it Weierstrass formula for surfaces of 
	prescribed mean curvature}, Math. Ann., {\bf 245}, (1979),89-99.
  \item \label{r10}
	Konopelchenko B.G., {\it Multidimensional integrable systems and 
	dynamics of surfaces in space}, preprint of 
	Institute of Mathematics, Taipei, August 1993; in: {\it National
	Workshop on Nonlinear Dynamics}, (Costato M., Degasperis A. and 
	Milani M., Eds.), Ital. Phys. Society, Bologna, 1995, pp. 33-40.
  \item \label{r11}
	Konopelchenko B.G., {\it Induced surfaces and their 
	integrable dynamics}, Stud. Appl. Math., {\bf 96}, (1996), 9-51.
  \item \label{r12}
        T.J. Willmore, {\it Total curvature in Riemannian Geometry},
        Ellis Horwood, New York, 1982.	
  \item \label{r13}
	Carroll R. and Konopelchenko B.G., {\it Generalized 
	Weierstrass-Enneper inducing, conformal immersions and gravity}, 
	Int. J. Modern Physics {\bf A11}, (1996), 1183-1216.
  \item \label{r14}
	Konopelchenko B.G. and Taimanov I., {\it Generalized Weierstrass 
	formulae, soliton equations and Willmore surfaces}, 
	preprint N. 187, Univ. Bochum, (1995).
  \item \label{r15}
        Konopelchenko B.G. and Taimanov I., {\it Constant mean 
	curvature surfaces via an integrable dynamical system}, 
	J. Phys. {\bf A}:Math. Gen., {\bf 29}, (1996), 1261-1265.
 \item \label{r16}
	Taimanov I., {\it Modified Novikov-Veselov equation and 
	differential geometry of surfaces}, Trans. Amer. Math. Soc., Ser. 2
	, {\bf 179}, (1997), 133-159.
 \item \label{r17}
	Grinevich P.G. and Schmidt M.V., {\it Conformal invariant 
	functionals of immersion of tori into ${\RR}^3$}, 
	Journal of Geometry and Physics (to appear); 
	preprint SFB288 N 291, TU-Berlin, 1997.
 \item \label{r18}
	Richter J., {\it Conformal maps of a Riemann surface into space of 
	quaternions}, PH.D Thesis, TU-Berlin, 1997.
 \item \label{r19}
        Taimanov I., {\it Global Weierstrass representation and its 
	spectrum}, Uspechi Mat. Nauk, {\bf 52}, N 6, (1997), 187-188.
 \item \label{r20}
        Taimanov I., {\it The Weierstrass representation of spheres in 
	${\RR}^3$, 
	the Willmore numbers and soliton spheres}, preprint SFB 288, 
	N 302, TU-Berlin, 1998.
 \item \label{r21}
	Matsutani S., {\it Immersion anomaly of Dirac operator for
	surfaces in $\RR^3$}, preprint, physics/9707003, (1997).
 \item \label{r22}
	Matsutani S., {\it On density state of quantized Willmore 
	surfaces-a way to quantized extrinsic string in $\RR^3$}, J.
	Phys. A: Math. Gen., {\bf 31}, (1998), 3595-3606.
 \item \label{r23}
        Konopelchenko B.G. and Landolfi G., {\it On classical string
	configurations}, Mod. Phys. Lett. A, {\bf 12}, (1997), 3161-3168.

 \item \label{r24}
	Bryant R.L., {\it Conformal and minimal immersions of compact
	surfaces into the $4-$sphere}, J. Diff. Geom., {\bf 17}, (1982)
	, 455-473.

\item \label{r25}
	Friedrich T., {\it On surfaces in four-spaces}, Ann. Glob. Analysis 
	and Geometry, {\bf 2}, (1984), 257-287.
\item \label{r26}
	Hoffman D.A. and Osserman R., {\it The Gauss map of surfaces in
	$\RR^3$ and $\RR^4$}, Proc. London Math. Soc. (3), {\bf 50}, 
	(1985), 27-56. 
\item \label{r27}
	Ferus D., Pedit F., Pinkall U. and Sterling I., {\it Minimal
	tori in $\SS^4$}, J. Reine Angew. Math., {\bf 429}, (1992), 1-47.
\item \label{r28}
	Baird P. and Wood J.C, {\it Weierstrass representations for
	harmonic morphisms on Euclidean spaces and surfaces}, preprint dg-ga
	/9512010, (1995).
\item \label{r29}
	Konopelchenko B.G. and Landolfi G., {\it Generalized Weierstrass 
	representation for surfaces in multidimensional Riemann spaces}, 
	math.DG/ 9804144, (1998); J. Geometry and Physics (to appear). 
\item \label{r30}
	Pedit F. and Pinkall U., in preparation.
\item \label{r31}
	Matsutani S., {\it Dirac operators of a conformal surface immersed
	in $\RR^4$: further generalized Weierstrass relation}, 
	solv-int/9801006, (1998).

\item	\label{r32}
	Barrett J., Gibbons G.W., Perry M.J. and Ruback P.,
	{\it Kleinian Geometry and the $N=2$ superstring}, Int. J.
	Modern Phys. A, {\bf 9}, (1994), 1457-1493.

 \item \label{r33}
	Bobenko A.I., {\it Surfaces in terms of $2$ by $2$ matrices. Old 
	and new integrable cases}, in: {\it Harmonic maps and integrable 
	systems}, (Fordy A. and Wood J., Eds.), pp. 83-127, Vieweg, 1994.
 \item \label{r34}
	Kusner R. and Schmitt N., {\it The spinor representation of surfaces 
	in space}, preprint dg-ga/9610005, 1996.

\item \label{r35}
	Kamberov G., Pedit F. and Pinkall U., {\it Bonnet pairs and 
	isothermic surfaces}, Duke Math. J., {\bf 92}, (1998), 637-644.
\item \label{r36}
	Ablowitz M.J. and Clarkson P.A., {\it Solitons, nonlinear evolution
	equations and inverse scattering}, Cambridge Univ. Press, Cambridge,
	1991.
\item \label{r37}
	Konopelchenko B.G., {\it Introduction to multidimensional 
	integrable equations}, Plenum Press, New York, 1992.
\item \label{r38}
	Konopelchenko B.G., {\it Solitons in multidimension}, World 
	Scientific, Singapore, 1993.
\item \label{r39}
	Hasimoto H., {\it A soliton on a vortex filament}, J. Fluid. Mech.,
	{\bf 51}, (1972), 477-485.

\end{enumerate}
\end{document}